\documentclass{amsart}

\usepackage{xypic}
\input xy
\xyoption{all}
\usepackage{epsfig}
\usepackage{amsthm}
\usepackage{amssymb}
\usepackage{amsmath}
\usepackage{amscd}

\newcommand{\bg}{\begin{equation}}
\newcommand{\ed}{\end{equation}}
\newcommand{\bga}{\begin{eqnarray}}
\newcommand{\eda}{\end{eqnarray}}

\newcommand{\lpf}{\textbf{Proof of Lemma:\ }}

\def\cbdu{\par{\raggedleft$\Box$\par}}

\newtheorem {Theorem}  {Theorem}

\numberwithin{Theorem}{section}

\newtheorem {Lemma}[Theorem]  {Lemma}

\theoremstyle{definition}
\newtheorem{Definition}[Theorem]{Definition}
\theoremstyle{remark}
\newtheorem{Remark}[Theorem]{Remark}

\expandafter\chardef\csname pre amssym.def
at\endcsname=\the\catcode`\@ \catcode`\@=11
\def\undefine#1{\let#1\undefined}
\def\newsymbol#1#2#3#4#5{\let\next@\relax
 \ifnum#2=\@ne\let\next@\msafam@\else
 \ifnum#2=\tw@\let\next@\msbfam@\fi\fi
 \mathchardef#1="#3\next@#4#5}
\def\mathhexbox@#1#2#3{\relax
 \ifmmode\mathpalette{}{\m@th\mathchar"#1#2#3}%
 \else\leavevmode\hbox{$\m@th\mathchar"#1#2#3$}\fi}
\def\hexnumber@#1{\ifcase#1 0\or 1\or 2\or 3\or 4\or 5\or 6\or 7\or 8\or
 9\or A\or B\or C\or D\or E\or F\fi}

\font\teneufm=eufm10 \font\seveneufm=eufm7 \font\fiveeufm=eufm5
\newfam\eufmfam
\textfont\eufmfam=\teneufm \scriptfont\eufmfam=\seveneufm
\scriptscriptfont\eufmfam=\fiveeufm

\catcode`\@=\csname pre amssym.def at\endcsname

\newcounter{remark}
\def  \12  {{\frac{1}{2}}}

\def\build#1_#2^#3{\mathrel{\mathop{\kern 0pt#1}\limits_{#2}^{#3}}}

\begin{document}

\title[Norm inflation for incompressible magneto-hydrodynamic system ]{Norm inflation for incompressible magneto-hydrodynamic system in $\dot{B}_{\infty}^{-1,\infty}$}


\author [Mimi Dai]{Mimi Dai}
\address{Department of Mathematics, UC Santa Cruz, Santa Cruz, CA 95064,USA}
\author[Jie Qing]{ Jie Qing}
\address{Department of Mathematics, UC Santa Cruz, Santa Cruz, CA 95064,USA}
\author[Maria E. Schonbek]{ Maria E. Schonbek}
\address{Department of Mathematics, UC Santa Cruz, Santa Cruz, CA 95064,USA}
\email{mdai@ucsc.edu} \email{qing@ucsc.edu}
\email{schonbek@ucsc.edu}

\thanks{The work of Mimi Dai was partially supported by NSF Grant
DMS-0700535 and DMS-0600692;}
\thanks{The work of Jie Qing was
partially supported by NSF Grant DMS-0700535;}
\thanks{The work of
M. Schonbek was partially supported by NSF Grant DMS-0600692}





\begin{abstract}
Based on the construction of Bourgain and Pavlovi\'{c} \cite{BP}, we
demonstrate that the solutions to the Cauchy problem for the three
dimensional incompressible magneto-hydrodynamics (MHD) system can
develop diferent types of norm inflations in $\dot{B}_{\infty}^{-1, \infty}$.
Particularly the magnetic field can develop norm inflation in short
time even when the velocity remains small and vice verse.
Efforts are made to present a very expository  development  of  the inginious construction
of Bourgain and Pavlovi\'{c} in \cite{BP}.

\bigskip


\hspace{0.02cm}CLASSIFICATION CODE: 76D03,35Q35.
\end{abstract}

\maketitle

\section{Introduction}

In this paper we consider the three dimensional incompressible
magneto-hydro -dynamics (MHD) system:
\begin{eqnarray}\label{MHD:PDE1}
&u_t- \triangle u +u\cdot\nabla u-b\cdot\nabla b+\nabla p=0,\\
&b_t- \triangle b+u\cdot\nabla b-b\cdot\nabla u=0,\notag\\
&\nabla \cdot u=0, \ \ \ \nabla \cdot b=0, \notag
\end{eqnarray}
with the initial conditions
\begin{equation}
u(x,0)=u_0(x),\;\;\; b(x,0)=b_0(x),
\end{equation}



where $x\in\mathbb{R}^3$, $t\geq 0$, $u$ is the fluid velocity, $b$
is the magnetic field. The initial data $u_0$ and $b_0$ are
divergence free. When the magnetic field $b(x,t)$ vanishes,
incompressible MHD system is just incompressible Navier-Stokes
equations. The solutions to MHD system also share the same scaling
property of solutions to Navier-Stokes equations, that is,
\begin{equation}
u_\lambda(x,t) =\lambda u(\lambda x, \lambda^2t), \ b_\lambda(x,t)
=\lambda b(\lambda x, \lambda^2t), \ p_\lambda(x,t) =\lambda^2
p(\lambda x, \lambda^2t) \notag
\end{equation}
solve the MHD system (\ref{MHD:PDE1}) with initial data
\begin{equation}
u_{0\lambda} =\lambda u_0(\lambda x), \ b_{0\lambda} =\lambda
b_0(\lambda x), \notag
\end{equation}
if $u(x, t)$ and $b(x, t)$ solve the MHD system (\ref{MHD:PDE1})
with the initial data $u_0(x)$ and $b_0(x)$. The spaces that are
invariant under the above scaling are called the critical spaces.
Examples of critical spaces in three dimension are
$$
\dot{H}^\frac 12 \hookrightarrow L^3 \hookrightarrow BMO^{-1}
\hookrightarrow \dot{B}^{-1, \infty}_\infty.
$$
(see \cite{Ca}, for example, for the discussions of the embeddings).

The study of Navier-Stokes equations as well as of MHD system in
critical spaces has been one of the focus of the research activities
since the initial work of Kato \cite{Kato}. In case of Navier-Stokes
equations, on one hand, in 2001, Koch and Tataru \cite{KT} were able
to establish  global well-posedness of Navier-Stokes equations with
small initial data in the space $BMO^{-1}$, on the other hand,
recently, Bourgain and Pavlovi\'{c} \cite{BP} showed the
ill-posedness for Navier-Stokes equations in $\dot{B}^{-1,\infty}_\infty$. More precisely, Bourgain and Pavlovi\'{c}
constructed some arbitrarily small initial data in $\dot{B}^{-1,\infty}_\infty$ and produced the so-called norm inflation in the sense that the solution becomes arbitrarily large in $\dot{B}^{-1,\infty}_\infty$ after an arbitrarily short time.

In a recent work \cite{MYZ}, Miao, Yuan and Zhang proved the
existence of a global mild solution in $BMO^{-1}$ for small initial
data and uniqueness of such solution in $C([0,\infty); BMO^{-1})$.
It is then an interesting problem to study the solutions to MHD
system with initial data in the space $\dot{B}^{-1, \infty}_\infty$.
In this paper, we  discuss different cases of norm inflation phenomena for
the MHD system in $\dot{B}_\infty^{-1,\infty}$. We  construct  arbitrarily small
initial data $(u_0, b_0)$  in
$\dot{B}_{\infty}^{-1, \infty}\times\dot{B}_{\infty}^{-1, \infty}$.  This data  when evolved in time  through the MHD system give raise to "norm inflation" in
$\dot{B}^{-1, \infty}_\infty$ for the  corresponding  solutions $(u, b)$. One particularly interesting scenario
is that the magnetic field $b$ shows norm inflation while the
velocity $u$ remains small. Namely, we show that
\begin{Theorem}\label{theorem}
For any $\delta >0$ there exists a solution $(u, b, p)$ to the MHD
system (\ref{MHD:PDE1}) with data $u_0\in\mathcal{S}$ and
$b_0\in\mathcal{S}$ which satisfy
\begin{equation}\notag
\|u(0)\|_{\dot{B}_{\infty}^{-1, \infty}} \lesssim \delta, \ \
\|b(0)\|_{\dot{B}_{\infty}^{-1, \infty}} \lesssim \delta,
\end{equation}
such that for some $0<T<\delta$
\begin{equation}\notag
\|b(T)\|_{\dot{B}_{\infty}^{-1, \infty}}\gtrsim \frac{1}{\delta}
\end{equation}
but for any $0<t<T<\delta$
\begin{equation}
\notag \|u(t)\|_{B_{\infty}^{-1, \infty}} \lesssim \delta.
\end{equation}
\end{Theorem}
\begin{Remark} We refer the reader to the  beginning of section of Preliminaries for the definition of the symbol $\lesssim$.
\end{Remark}
 Our proof follows the methods  introduced by Bourgain and
Pavlovi\'c in \cite{BP}. 
Efforts are made to give a very expository  development  of 
 the inginious ideas of Bourgain and
Pavlovi\'{c} in \cite{BP}.

We now recall some  auxiliary concepts  related to plane waves, which are necessary in the sequel:
\begin{itemize}
 \item  The ``diffusion" of a plane wave $v\sin (k\cdot x)$ in $R^3$ is given by
\[e^{\Delta t}v\sin (k\cdot x)=e^{-|k|^2 t}v\sin (k\cdot x)\] Thus the magnitude of the diffusion of a
plane wave dies down in time in the scale that is measured by the
square of the magnitude of the wave vector $k$.
\item It is easy to  see
that $u = b = e^{-|k|^2 t}v\sin (k\cdot x)$ solve MHD system when the
wave vector $k$ is orthogonal to the amplitude vector $v$.
\item The
nonlinear interaction of two such diffusions in MHD system can be
captured, and it only produces a slower diffusion if the two wave
vectors are close.
\end{itemize}
We note that these observations, are the basis of the original argument  of Bourgain and
Pavlovi\'{c}  in \cite{BP}. We will use them to  construct a
combination of such ``diffusions"  with least nonlinear interactions yet
producing enough slower  ``diffusions"  to cause the norm inflation in
short time.
\begin{Remark} \label{remark1}
It is interesting to observe that even though the initial velocity
is zero the velocity can be triggered to develop norm inflation
while the magnetic field stays under control. More precisely we can
show that, for any $\delta >0$ there exists a solution $(u, b, p)$
to the MHD system (\ref{MHD:PDE1}) with vanishing initial velocity
and some $b_0\in\mathcal{S}$ which satisfies that
$$
\|b(0)\|_{\dot{B}_{\infty}^{-1, \infty}} \lesssim \delta
$$
and that for  some $0<T<\delta$
\begin{equation}\notag
\|u(T)\|_{\dot{B}_{\infty}^{-1, \infty}} \gtrsim 1/ \delta,
\end{equation}
while for all $0<t<T<\delta$
\begin{equation}\notag
\|b(t)\|_{B_{\infty}^{-1, \infty}} \lesssim \delta.
\end{equation}
\end{Remark}
\begin{Remark}\label{remark2}
We also note  that due to the interaction
between the velocity and the magnetic field,  if initially they are the same, they may restrain each
other from norm inflations.
\end{Remark}
In our paper, we  present our results in $\mathbb{T}^3$. But, as
pointed out in \cite{BP}, the proof can be modified to
$\mathbb{R}^3$.

\bigskip

\section{Preliminaries}

\subsection{Notation}
We denote by $A\lesssim B$ an estimate of the form $A\leq C B$ with
some constant $C$, and by $A\sim B$ an estimate of the form $C_1
B\leq A\leq C_2 B$ with some constants $C_1$, $C_2$. For
completeness we recall the  defining norms for the Besov space
$\dot{B}_{\infty}^{-1, \infty}$ and the $BMO^{-1}$ space
\begin{equation}\label{hombes}
\|f\|_{\dot{B}_{\infty}^{-1, \infty}}=\sup_{t>0} t^{1/2}\|
e^{t\triangle}f\|_{L^\infty}.
\end{equation}
\begin{equation}
\|f\|_{BMO^{-1}}=\sup_{x_0\in R^3, R>0}\left(\frac{1}{|B(x_0,
\sqrt{R})|} \int\limits_{0}^{R}\int\limits_{B(x_0, \sqrt{R})}
|e^{t\triangle}f(y)|^2\, dy\, dt\right)^{\frac{1}{2}}.
\end{equation}
We will also work with the so-called inhomogeneous Besov space
$B^{-1, \infty}_\infty$ with the norm
\begin{equation}\label{inhombes}
\|f\|_{B^{-1, \infty}_\infty} = \sup_{0<t<1}\sqrt t
\|e^{t\Delta}f\|_{L^\infty}.
\end{equation}
Clearly
$$
\|f\|_{B^{-1, \infty}_\infty} \leq \|f\|_{\dot{B}_{\infty}^{-1,
\infty}}.
$$
and,
$$
\|f\|_{B^{-1, \infty}_\infty} \leq \|f\|_{L^\infty},
$$
since $\|e^{t\Delta}f\|_{L^\infty} \leq \|f\|_{L^\infty}$.

\subsection{The well-posedness result of the incompressible MHD system in $BMO^{-1}$}
We recall the well-posedness result of C. Miao, B. Yuan and B. Zhang
in $BMO^{-1}$ in \cite{MYZ}. For this we introduce the spaces $X_T$
and the corresponding norm.

\begin{Definition}
Let $u(x, t)$ be a measurable function on $\mathbb{R}^3\times [0,
T)$ for $T>0$ and let
\begin{eqnarray}\label{XTnorm}
& \|u(\cdot, \cdot)\|_{X_T}  = \sup_{0<t<T} t^{1/2}\|u(\cdot,
t)\|_{L^\infty} \\ & + \sup_{x_0\in R^3,
0<R<T}\left(\frac{1}{|B(x_0,
\sqrt{R})|}\int\limits_{0}^{R}\int\limits_{B(x_0, \sqrt{R})} |u(y,
t)|^2\, dy\, dt\right)^{\frac{1}{2}}.\notag
\end{eqnarray}
Then the space-time space $X_T$ is defined by
\begin{equation}
X_T = \{f(x, t)\in L^2(0, T; L^2(\mathbb{R}^3)):
\|f\|_{X_T}<\infty\}\notag
\end{equation}
\end{Definition}
It is worth to mention that, for each $t\in (0, T]$,
$$
\|f(\cdot, t)\|_{L^\infty} \leq \frac 1{\sqrt t} \|f\|_{X_T}.
$$
In \cite{MYZ}, Miao, Yuan and Zhang  proved the following existence
theorem:
\begin{Theorem} (Miao, Yuan and Zhang)
Let $(u_0(x), b_0(x))\in BMO^{-1}\times BMO^{-1}$ with $\nabla\cdot u_0=0\; \mbox{and}\;
\nabla\cdot b_0=0$. Then, there exists a positive constant
$\varepsilon $ such that if $\|(u_0, b_0)\|_{BMO^{-1}}<\varepsilon $
then the MHD system has a unique global mild solution $(u(x, t),b(x, t))\in X_T\times X_T$ satisfying that $\|(u(x, t), b(x, t))\|_{X_T\times X_T} \leq 2\varepsilon$ for all $T>0$.
\end{Theorem}

\subsection{Bilinear operators}
Let $\mathbb{P}$ denote the projection on divergence-free vector
fields, which acts on a function $\phi$ as
\begin{equation}\notag
\mathbb{P}(\phi)=\phi+\nabla\cdot(-\triangle)^{-1}div\phi.
\end{equation}
As shown in \cite{KT,MYZ} the bilinear operator
\begin{equation}\notag
\mathcal{B} (u, v)=\int\limits_{0}^{t} e^{(t-\tau)\triangle}
\mathbb{P}(u\cdot \nabla v)\, d\tau ,
\end{equation}
maps $X_T\times X_T$ into $X_T$ continuously. More precisely we
have,
\begin{equation}\label{B}
\|\mathcal{B} (u, v)\|_{X_T}\lesssim \|u\|_{X_T}\|v\|_{X_T}.
\end{equation}

\subsection{Rewriting the MHD system}
Following ideas from \cite{BP} we rewrite the MHD system
(\ref{MHD:PDE1})  introducing the  expression
\begin{equation}\label{u}
u=e^{t\triangle}u_0-u_1+y
\end{equation}
\begin{equation}\label{b}
b=e^{t\triangle}b_0-b_1+z
\end{equation}
where
\begin{equation} \label{u1}
u_1(x, t)= \mathcal{B}(e^{t\triangle}u_0(x),e^{t\triangle}u_0(x))
-\mathcal{B}(e^{t\triangle}b_0(x),e^{t\triangle}b_0(x)),
\end{equation}
\begin{equation} \label{b1}
b_1(x, t)= \mathcal{B}(e^{t\triangle}u_0(x),e^{t\triangle}b_0(x))
-\mathcal{B}(e^{t\triangle}b_0(x),e^{t\triangle}u_0(x)),
\end{equation}
An easy calculation shows that
\begin{align} \label{yz}
y_t-\triangle y+G_0 + G_1 + G_2=0,\\
z_t-\triangle z+K_0 + K_1 + K_2=0,\notag
\end{align}
where
\begin{eqnarray}
&G_0=\mathbb{P}[(e^{t \triangle} u_0\cdot\nabla)u_1+(u_1\cdot\nabla)e^{t \triangle} u_0+(u_1\cdot\nabla)u_1]\notag\\
&-\mathbb{P}[(e^{t \triangle}
b_0\cdot\nabla)b_1+(b_1\cdot\nabla)e^{t \triangle}
b_0+(b_1\cdot\nabla)b_1]\notag
\end{eqnarray}
\begin{eqnarray}
&G_1=\mathbb{P}[(e^{t \triangle} u_0\cdot\nabla)y+(u_1\cdot\nabla)y+(y\cdot\nabla)e^{t \triangle} u_0+(y\cdot\nabla)u_1]\notag\\
&-\mathbb{P}[(e^{t \triangle}
b_0\cdot\nabla)z+(b_1\cdot\nabla)z+(z\cdot\nabla)e^{t \triangle}
b_0+(z\cdot\nabla)b_1]\notag
\end{eqnarray}
\begin{equation}
G_2=\mathbb{P}[(y\cdot\nabla)y]-\mathbb{P}[(z\cdot\nabla)z]\notag
\end{equation}
and
\begin{eqnarray}
&K_0=\mathbb{P}[(e^{t \triangle} u_0\cdot\nabla)b_1+(u_1\cdot\nabla)e^{t \triangle} b_0+(u_1\cdot\nabla)b_1]\notag\\
&-\mathbb{P}[(e^{t \triangle}
b_0\cdot\nabla)u_1+(b_1\cdot\nabla)e^{t \triangle}
u_0+(b_1\cdot\nabla)u_1]\notag
\end{eqnarray}
\begin{eqnarray}
&K_1=\mathbb{P}[(e^{t \triangle} u_0\cdot\nabla)z+(u_1\cdot\nabla)z+(y\cdot\nabla)e^{t \triangle} b_0+(y\cdot\nabla)b_1]\notag\\
&-\mathbb{P}[(e^{t \triangle}
b_0\cdot\nabla)y+(b_1\cdot\nabla)y+(z\cdot\nabla)e^{t \triangle}
u_0+(z\cdot\nabla)u_1]\notag
\end{eqnarray}
\begin{equation}
K_2=\mathbb{P}[(y\cdot\nabla)z]-\mathbb{P}[(z\cdot\nabla)y].\notag
\end{equation}
Here $G_0$ and $K_0$ are constants, $G_1$ and $K_1$ are linear,
while $G_2$ and $K_2$ are quadratic in terms of $y$ and $z$.
\begin{Remark}
Note that  although the second equation in the MHD system (\ref{MHD:PDE1}) has no pressure, since $u$ and $b$ are both divergence free,  the term $u\cdot\nabla b-b\cdot\nabla u$ is automatically divergence free. Hence the projector $\mathbb{P}$ acting on this term does not change the second equation and hence we can write $b_1$ and $K_i$'s as  described above. 

\end{Remark}

\bigskip

\section{Interactions of plane waves}
In this section we show how the diffusions of plane
waves interact in MHD system.  These interactions are the basis for the
constructions of initial data which will  evolve into  the different  cases  of velocity and  magnetic field
 norm inflations .
\subsection{Diffusion of a plane wave } As a first step,
we consider the same initial data for velocity and magnetic using
one single plane wave. Suppose $k\in \mathbb{R}^3$, $v \in\mathbb{S}^2$ and $k \cdot v=0$. Let
$$
u_0 = b_0 = v\cos(k\cdot x).
$$
Then
$ \nabla\cdot u_0 = 0,\,\mbox{ and}\;\nabla\cdot b_0 = 0. $ and
\begin{equation}\label{difu}
e^{t\triangle}v\cos(k\cdot x)=e^{-|k|^2t}v\cos(k\cdot x).
\end{equation}

In fact the ``diffusions"  $(e^{t\Delta}v\cos(k\cdot x),e^{t\Delta}v\cos(k\cdot x))$ of a plane wave solve the MHD system with vanishing pressure. It is important to notice that
\begin{itemize}
\item $
\|v\cos(k\cdot x)\|_{\dot{B}^{-1, \infty}_\infty} \sim \frac 1{|k|},
$
\item$
\|e^{t\Delta}v\cos (k\cdot x)\|_{X_T} \lesssim \frac 1{|k|}.
$
\end{itemize}
which says that the size of a plane wave in the space $\dot{B}^{-1,\infty}_\infty$  is reciprocal to the magnitude of its wave vector, and in $X_T$ it is bounded by this same reciprocal.

\subsection{Interaction of plane waves}

Now we consider  different plane wave initial data for the velocity and magnetic. Suppose $k_i\in\mathbb{R}^3$, $v_i \in\mathbb{S}^2$ and $k_i\cdot v_i=0$, for $i=1, 2$. Let
\begin{equation}\notag
u_0= \cos(k_1\cdot x)v_1,
\end{equation}
\begin{equation}\notag
b_0= \cos(k_2\cdot x)v_2.
\end{equation}
Using the decomposition given in Section 2.4,
$$
u = e^{t\Delta}u_0, \quad b = e^{t\Delta}b_0 - b_1 + z
$$
solve the MHD system with vanishing pressure. To simplify our
calculations we  assume that
$$
k_2\cdot v_1 = 0, \ \ \text{and} \ \  k_1\cdot v_2 = \frac 12,
$$
which eliminates the term $e^{t\Delta}u_0\cdot \nabla (e^{t\Delta}b_0)$ and
gives
\begin{align}\notag
e^{t \triangle} b_0\cdot\nabla(e^{t \triangle} u_0)
&=- e^{-(|k_1|^2+|k_2|^2) t}v_1\sin(k_1\cdot x)\cos(k_2\cdot x)(k_1\cdot v_2)\\
&=-\frac{1}{4} e^{-(|k_1|^2+|k_2|^2)t}v_1 (\sin((k_1-k_2)\cdot x) +
\sin((k_1+k_2)\cdot x)). \notag
\end{align}
Hence
$$
\aligned b_1 & = \frac 14  v_1\sin((k_1-k_2)\cdot x) \int_0^t
e^{-(|k_1|^2+|k_2|^2)\tau}e^{-|k_1-k_2|^2(t-\tau)}d\tau \\ & + \frac
14 v_1 \sin((k_1+k_2)\cdot x)\int_0^t
e^{-(|k_1|^2+|k_2|^2)\tau}e^{-|k_1+k_2|^2(t-\tau)}d\tau \\ & = b_{1,
0} + b_{1, 1},\endaligned
$$
where
$$
b_{1, 0} = \frac 14  v_1 \sin((k_1-k_2)\cdot
x)\frac{-e^{-(|k_1|^2+|k_2|^2)t}+e^{-|k_1-k_2|^2t}}{2k_1\cdot k_2}
$$
and
$$
b_{1,1} =  \frac 14 v_1 \sin((k_1+k_2)\cdot
x)\frac{e^{-(|k_1|^2+|k_2|^2)t}-e^{-|k_1+k_2|^2t}}{2k_1\cdot k_2}.
$$
Therefore, if we can manage to control $z$ in the light of the
continuity of the bilinear operator $\mathcal{B}$ in $X_T$, then the
interaction of two plane waves are small in
$\dot{B}^{-1,\infty}_\infty$ if neither the sum nor the difference
of their wave vectors is small in magnitude. In the mean time, the
interaction is sizable in $\dot{B}^{-1,\infty}_\infty$ if either the
sum or the difference of their wave vectors is small in magnitude.

\bigskip

\section{Proof of theorem \ref{theorem}}
In this section we will follow the idea from \cite{BP} to construct
initial data to produce norm inflation for solutions to MHD system.
From the discussions in the previous sections we know that the
interaction of two plane waves is not enough to show the norm
inflation. We need to build interactions of more plane waves. The
construction in \cite{BP} depends on the rather sophisticated
choices of plane waves. We will use a similar scheme.

\subsection{Construction of initial data for the MHD system}
For a fixed small number $\delta>0$ we will specify later  following initial data
:
\begin{equation}\label{u0}
u_0=\frac{Q}{\sqrt{r}}\sum_{s=1}^r|k_s|v_s \cos(k_s\cdot x)
\end{equation}
and
\begin{equation}\label{b0}
b_0=\frac{Q}{\sqrt{r}}\sum_{s=1}^r|k_s'|v_s' \cos(k_s'\cdot x).
\end{equation}
We expect for each $s$, the interaction of the two plane waves
$v_s\cos(k_s\cdot x)$ and $v_s'\cos(k_s'\cdot x)$ is sizable in
$\dot{B}^{-1,\infty}_\infty$; while the interactions of plane waves
of different $s$ is small, as demonstrated in Section 3.2. Hence
\begin{itemize}
\item{Wave vectors:} The wave vectors $k_s\in\mathbb{R}^3$ are parallel to a given vector
$k_0\in\mathbb{R}^3$. The modulo $|k_0|$ will be taken to be large,
depending on $Q$. The magnitude of $k_s$ is defined by,
\begin{align}\label{kl}
|k_s|=2^s|k_0||k_{s-1}|, \ \ \ \ s=1, 2, 3, ..., r.
\end{align}
The wave vectors $k'_s\in\mathbb{R}^3$ is defined by
\begin{equation} \label{ks-ks'}
k_s-k'_s=\eta
\end{equation}
for a given vector $\eta\in \mathbb{S}^2$.
\item{Amplitude vectors:} The amplitude vectors $v_s, v'_s\in\mathbb{S}^2$ satisfy
\begin{equation} \label{v-equation}
k_s\cdot v_s=k'_s\cdot v'_s=0
\end{equation}
to ensure the initial data are divergence free.
\item{Auxiliary assumptions:} We also require that
\begin{equation} \label{eta-v}
\eta\cdot v_s=0, \ \ \ \eta\cdot v'_s=\frac{1}{2}
\end{equation}
to simplify our calculations. In fact we will choose $v_s =v$ a
fixed vector.
\end{itemize}
We first point out  the following  simple
facts to further motivate the choices of the magnitudes of $k_s$.
\begin{Lemma}\label{k}
\begin{equation}
\sum_{l<s}|k_l|\sim |k_{s-1}| \quad\text{and} \ \sum_{l<s}|k_l'|\sim
|k_{s-1}'|
\end{equation}
\begin{equation}
\sum_{s=1}^r |k_s|e^{-|k_s|^2t} \lesssim \frac 1{\sqrt t} \quad
\text{and} \ \sum_{s=1}^r |k_s'|e^{-|k_s'|^2t} \lesssim \frac
1{\sqrt t}
\end{equation}
\begin{equation}
v_i \cdot k_j = v_i\cdot k_j' = v_i\cdot \eta = 0, \quad \forall
\quad i,j = 1,2, \dots, r.
\end{equation}
\begin{equation}
\sum_{i=1}^r |k_i|e^{-\frac {|k_i|^2}{|k_0|^2}} \lesssim 1,
\quad\text{and} \ \sum_{i=1}^r |k_i'|e^{-\frac {|k_i'|^2}{|k_0|^2}}
\lesssim 1.
\end{equation}
\end{Lemma}
\lpf By the definition (\ref{kl}), it is clear that
$|k_{l-1}|<\frac{1}{2}|k_l|$, which easily implies the first
statement. For second statement, again due to the definition
(\ref{kl}), we know that $|k_s|\sim |k_s|-|k_{s-1}|$. Thus,
$$
\sum_{s=1}^r|k_s|e^{-|k_s|^2t} \sim
\sum_{s=1}^r(|k_s|-|k_{s-1}|)e^{-|k_s|^2t},
$$
while the later one can be considered as a finite Riemman summation
of the function $e^{-x^2t}$. Therefore
$$
\sum_{s=1}^r|k_s|e^{-|k_s|^2t} \lesssim \int_0^\infty e^{-x^2t} dx =
\frac{1}{\sqrt{t}} \int_0^\infty e^{-x^2t}d(x\sqrt{t}) = \frac{\sqrt
\pi}{2\sqrt{t}}.
$$
\noindent For the third statement, we note that, $k_s$ are parallel to the given vector $k_0$ for all $s$ and $v_s=v$ is a fixed vector by the above choice. Hence, by (\ref{v-equation}) 
$$
v_i\cdot k_j=0, \ \ \mbox { for all } i,j=1,2, ..., r. 
$$
On the other hand side,  from (\ref{ks-ks'}) and (\ref{eta-v}), we have
$$
v_i\cdot k'_j=v_i\cdot (k_j-\eta)=v_i\cdot k_j-v_i\cdot\eta=0.
$$
The forth statement in the Lemma follows from the second one provided an appropriate choice of $k_0$. It completes the proof of Lemma.
\cbdu
Next we calculate the norm of our initial data.
\begin{Lemma} For $u_0$ and $b_0$ given in (\ref{u0}) (\ref{b0}) we have
\begin{equation}\label{norm:u0-1}
\|u_0\|_{\dot{B}_{\infty}^{-1, \infty}}\lesssim\frac{Q}{\sqrt{r}}, \
\ \|b_0\|_{\dot{B}_{\infty}^{-1, \infty}}\lesssim\frac{Q}{\sqrt{r}}.
\end{equation}
\end{Lemma}
\lpf For the given initial data $u_0$, we have that, due to(\ref{difu})
\begin{equation} \label{h-u0}
e^{\tau\triangle}u_0=\frac{Q}{\sqrt{r}}\sum_{s=1}^r|k_s|v_s
\cos(k_s\cdot x)e^{-|k_s|^2\tau}.
\end{equation}
Hence by Lemma \ref{k}
\begin{align} \label{first-bound}
\|u_0\|_{\dot{B}_{\infty}^{-1, \infty}} &\sim \frac{Q}{\sqrt{r}}
\sup_{t>0} \sqrt{t}\sum_{s=1}^r |k_s|e^{-|k_s|^2t}\notag\\
&\lesssim \frac{Q}{\sqrt{r}},
\end{align}
The bound for $\|b_0\|_{\dot{B}_{\infty}^{-1, \infty}}$ follows in
the similar way. \cbdu
\begin{Lemma}\label{XTu0} For $u_0$ and $b_0$ given in (\ref{u0}) (\ref{b0}) we have
\begin{equation}\label{norm:u0-2}
\|e^{t\Delta}u_0\|_{X_T} \lesssim Q, \ \ \|e^{t\Delta}b_0\|_{X_T}
\lesssim Q.
\end{equation}
\end{Lemma}
\lpf We only need to verify one of the two.
$$
\|e^{t\Delta} u_0\|_{X_T} \lesssim \frac {Q}{\sqrt r}\left(1 + \sup_{t\in
[0, T]}(\int_0^t (\sum_{i=1}^r |k_i|e^{-|k_i|^2\tau})^2 d\tau )^\frac{1}{2}\right),
$$
where
$$
(\sum_{i=1}^r |k_i|e^{-|k_i|^2\tau})^2 \lesssim \sum_{i=1}^r |k_i|^2
e^{-2|k_i|^2\tau} + 2\sum_{i=1}^r
|k_i|e^{-|k_i|^2\tau}\sum_{j<i}|k_j|\lesssim \sum_{i=1}^r |k_i|^2
e^{-|k_i|^2\tau}
$$
Hence 
$$
\int_0^t (\sum_{i=1}^r |k_i|^2 e^{-|k_i|^2 \tau})^2d\tau ) \lesssim
\sum_{i=1}^r (1 - e^{-|k_i|^2t})\lesssim r,
$$
which implies
$$
\|e^{t\Delta} u_0\|_{X_T} \lesssim  \frac Q{\sqrt r} + Q .
$$
This completes the proof of the Lemma. \cbdu
Finally we make a note that
\begin{Lemma} For $t\in [0, +\infty)$,
\begin{equation}\label{norm:u0-3}
\|e^{t\Delta}u_0\|_{\dot{B}^{-1, \infty}_\infty} \lesssim \frac
Q{\sqrt r}e^{-|k_0|^2t}, \ \ \|e^{t\Delta}b_0\|_{\dot{B}^{-1,
\infty}_\infty} \lesssim \frac Q{\sqrt r}e^{-|k_0|^2t}.
\end{equation}
\end{Lemma}

\bigskip

\subsection{Analysis of $u_1$}

As demonstrated in Section 3.2 we consider the decomposition
$$
\aligned u  & = e^{t\Delta}u_0 - u_1 + y \\
b & = e^{t\Delta}b_0 - b_1 + z\endaligned
$$
We want to handle the $u_1$ first. Recall the definition (\ref{u1})
$$
u_1 = \mathcal{B}(e^{t\Delta}u_0, e^{t\Delta}u_0) -
\mathcal{B}(e^{t\Delta}b_0, e^{t\Delta}b_0)
$$
By our discussions in Section 3.2 the interactions should be small.
By the fact that $v_i\cdot k_j=0$ it is immediately seen that
$$
e^{t\triangle}u_0\cdot\nabla e^{t\triangle}u_0  = 0.
$$
Then a straight calculation shows
$$
\aligned e^{t\triangle}b_0\cdot\nabla e^{t\triangle}b_0 & =
-\frac{Q^2}{r}\sum_{i,j=1}^r|k_i'||k_j'|e^{-(|k_i'|^2+|k_j'|^2)t}(v_i'\cdot
k_j')v_j'\cos(k_i'\cdot x)\sin(k_j'\cdot x)\\ & = -\frac {Q^2}{2r}
\sum_{i,j=1}^r|k_i'||k_j'|e^{-(|k_i'|^2+|k_j'|^2)t}(v_i'\cdot
k_j')v_j'\sin
(k_i'+k_j')\cdot x \\
& \quad -\frac {Q^2}{2r}
\sum_{i,j=1}^r|k_i'||k_j'|e^{-(|k_i'|^2+|k_j'|^2)t}(v_i'\cdot
k_j')v_j'\sin (k_j'-k_i')\cdot x\endaligned
$$
and
$$
\aligned \notag\mathbb{P}e^{t\triangle}b_0\cdot\nabla e^{t\triangle}b_0 &
= -\frac {Q^2}{2r}
\sum_{i,j=1}^r|k_i'||k_j'|e^{-(|k_i'|^2+|k_j'|^2)t}(v_i'\cdot
k_j')u_j\sin (k_j'+k_i')\cdot x\notag \\
& \quad -\frac {Q^2}{2r}
\sum_{i,j=1}^r|k_i'||k_j'|e^{-(|k_i'|^2+|k_j'|^2)t}(v_i'\cdot
k_j')w_j\sin (k_j'-k_i')\cdot x\\ & = E_1 + E_2,\notag\endaligned
$$
where $u_j$ is the projection of $v_j'$ to the orthogonal to
$k_j'+k_i'$ and $w_j$ is the projection of $v_j'$ to the orthogonal
to $k_j'-k_i'$. Hence
$$
\mathcal{B}((e^{t\Delta}b_0, e^{t\Delta}b_0)  = \int_0^t
e^{(t-\tau)\Delta}E_1d\tau + \int_0^t e^{(t-\tau)\Delta}E_2d\tau =
F_1 + F_2. $$
\noindent We work for
$F_1$ in the following. $F_2$ will be handled similarly and so the detail is omitted.
$$
F_1 = \frac {Q^2}{2r} \sum_{i,j=1}^r|k_i'||k_j'|(v_i'\cdot
k_j')u_j\sin (k_j'+k_i')\cdot x
\frac{e^{-(|k_i'|^2+|k_j'|^2)t}-e^{-|k_i'+k_j'|^2t}}{|k_i'+k_j'|^2
- (|k_i'|^2 + |k_j'|^2)},
$$
By the fact that $k_i'\cdot v_i' = 0$ and since the function $\frac {1 - e^{-x}}x$ is bounded for $x>0$, we have
$$
|e^{\tau\Delta}F_1| \lesssim \frac {Q^2}r \sum_{j=1}^r \sum_{i<j}
|k_i'||k_j'|^2te^{-(|k_i'|^2+|k_j'|^2)t} e^{-|k_i'+k_j'|^2\tau}
$$
and
$$
|F_1| \lesssim \frac {Q^2}r \sum_{j=1}^r \sum_{i<j}
|k_i'||k_j'|^2te^{-(|k_i'|^2+|k_j'|^2)t}.
$$
 Hence
$$
|F_1| \lesssim \frac {Q^2}r \sum_{j=1}^r \sum_{i<j}
|k_j'|^2te^{-\frac 14 |k_j'|^2t}|k_i'|e^{-\frac 14
|k_j'|^2t}\lesssim \frac {Q^2}r \sum_{j=1}^r |k_{j-1}'|e^{-\frac 14
|k_j'|^2t}.
$$
and
$$
|e^{\tau\Delta}F_1| \lesssim \frac {Q^2}r \sum_{j=1}^r
|k_{j-1}'|e^{-\frac 14 |k_j'|^2t}e^{-\frac 12|k_j'|^2\tau},
$$
where we use the fact that $xe^{-x}$ is bounded for $x>0$ and Lemma
\ref{k}.
$$
\aligned \|F_1\|_{X_T} & \lesssim \frac {Q^2}{r} \sum_{j=1}^r
|k_{j-1}'| \sup_{t\in [0, T]}\sqrt t e^{- \frac 14|k_j'|^2t} \\ &
\quad  + \sup_{t\in [0, T]} (\int_0^t (\sum_{j=1}^r
|k_{j-1}'|e^{-\frac 14 |k_j'|^2\tau})^2 d\tau)^\frac 12
\endaligned
$$
where
$$
\sum_{j=1}^r |k'_{j-1}|\sup_{t\in [0, T]} \sqrt t e^{- \frac
14|k_j'|^2t} \lesssim \sum_{j=1}^r \frac
{|k_{j-1}'|}{|k_j'|}\lesssim 1
$$
and
$$
\sup_{t\in [0, T]} (\int_0^t (\sum_{j=1}^r |k_{j-1}'|e^{-\frac 14
|k_j'|^2\tau})^2 d\tau)^\frac 12 \lesssim (\sum_{j=1}^r \frac
{|k_{j-1}'|}{|k_j'|})^\frac 12 \lesssim 1
$$
by the same argument as used in the proof of Lemma \ref{XTu0}. And
similarly
$$
\|F_1\|_{\dot{B}^{-1, \infty}_\infty} = \sup_{\tau > 0} \sqrt \tau
\|e^{\tau\Delta}F_1\|_{L^\infty} \lesssim \frac {Q^2}r \sum_{j=1}^r
\frac {|k_{j-1}'|}{|k_j'|} \lesssim \frac {Q^2}r.
$$
Therefore we conclude that
\begin{Lemma}\label{XYu1}
\begin{equation}\label{XTu1}
\|u_1\|_{\dot{B}^{-1, \infty}_\infty} \lesssim \frac {Q^2}r,
\quad \|u_1\|_{X_T} \lesssim \frac {Q^2}r.
\end{equation}
\end{Lemma}
\lpf $F_2$ can be handled just like what we did with $F_1$. \cbdu

\subsection{Analysis of $b_1$} First we recall that from (\ref{b1})
\begin{equation}\notag
b_1(x, t)= \mathcal{B}(e^{t\triangle}u_0(x),e^{t\triangle}b_0(x))
-\mathcal{B}(e^{t\triangle}b_0(x),e^{t\triangle}u_0(x)).
\end{equation}
Similar to the calculations in the previous section, first, due to
the fact that $v_i\cdot k_j' =0$, we have
$$
e^{\tau\triangle}u_0\cdot\nabla e^{\tau\triangle}b_0  = 0.
$$
And
$$
\aligned e^{\tau\triangle}b_0\cdot\nabla e^{\tau\triangle}u_0 & =
-\frac {Q^2}{2r}
\sum_{i,j=1}^r|k_i'||k_j|e^{-(|k_i'|^2+|k_j|^2)t}(v_i'\cdot
k_j)v_j\sin
(k_j+k_i')\cdot x \\
& \quad -\frac {Q^2}{2r} \sum_{i\neq
j}^r|k_i'||k_j|e^{-(|k_i'|^2+|k_j|^2)t}(v_i'\cdot
k_j)v_j\sin (k_j-k_i')\cdot x\\
& \quad -\frac {Q^2}{2r}
\sum_{i=1}^r|k_i||k_i'|e^{-(|k_i|^2+|k_i'|^2)t}(v_i'\cdot
k_i)v_i\sin (k_i-k_i')\cdot x.\endaligned
$$
We will write
$$
e^{\tau\triangle}u_0\cdot\nabla
e^{\tau\triangle}b_0-e^{\tau\triangle}b_0\cdot\nabla
e^{\tau\triangle}u_0= \tilde b_{1,0} + \tilde b_{1,1},
$$
where
$$
\aligned \tilde b_{1,0} & = \frac{Q^2}{2r} \sum_{i=1}^r|k_i||k_i'|
e^{-(|k_i|^2+|k_i'|^2)t} \sin(k_i-k_i')\cdot x
(v_i'\cdot k_i)v_i\\
& = \frac{Q^2}{4r} \sin (\eta\cdot x) \sum_{i=1}^r|k_i||k_i'|
e^{-(|k_i|^2+|k_i'|^2)t}v_i,\endaligned
$$
due to our choices of wave vectors and amplitude vectors in Section
4.1. We then set
$$
\tilde b_{1,1} = \tilde d_{1,1} + \tilde e_{1,1},
$$
where
$$
\aligned \tilde d_{1,1} & =  \frac{Q^2}{2r} \sum_{i=1}^r|k_i||k_i'|
e^{-(|k_i|^2+|k_i'|^2)t} \sin(k_i + k_i')\cdot x (v_i'\cdot k_i)v_i\\
& = \frac{Q^2}{4r} \sum_{i=1}^r|k_i||k_i'| e^{-(|k_i|^2+|k_i'|^2)t}
v_i \sin(k_i + k_i')\cdot x.\endaligned
$$
By the choices of $k_i'$, which behaves more or less like $k_i$ for
each $i$ when $|k_0|$ is very large, we conclude that $\tilde
e_{1,1}$ can be handled just like what we did for $E_1$ in the
previous section. We then have
$$
\aligned d_{1,1} & = \int_0^t e^{\Delta (t-\tau)}\mathbb{P}\tilde
d_{1,1}(\tau)d\tau \\
& = \frac {Q^2}{4r} \sum_{i=1}^r|k_i||k_i'|v_i \sin(k_i + k_i')\cdot
x \frac {e^{-(|k_i|^2+|k_i'|^2)t} - e^{-|k_i+k_i'|^2t}}{|k_i+k_i'|^2
- (|k_i|^2 + |k_i'|^2)}\endaligned
$$
which gives
$$
|d_{1,1}| \lesssim \frac {Q^2}r \sum_{i=1}^r |k_i|^2 te^{-\frac
12|k_i|^2t} \lesssim \frac {Q^2}r \sum_{i=1}^r e^{-\frac 14|k_i|^2t}
$$
and
$$
|e^{\tau\Delta} d_{1,1}|\lesssim \frac {Q^2}r \sum_{i=1}^r e^{-\frac
14|k_i|^2t} e^{-\frac 12 |k_i|^2 \tau}.
$$
Hence it is even easier to handle $d_{1,1}$ than to handle $F_1$ in
the previous section. Now, if denote
$$
b_{1,1} = \int_0^t e^{\Delta (t-\tau)}\mathbb{P}\tilde
b_{1,1}(\cdot, \tau)d\tau,
$$
we may conclude that
\begin{Lemma}
\begin{equation}\label{b:b11}
\|b_{1,1}\|_{\dot{B}^{-1, \infty}_\infty} \lesssim \frac {Q^2}r,
\quad \|b_{1,1}\|_{X_T} \lesssim \frac {Q^2}r.
\end{equation}
\end{Lemma}
The focus is now on
\begin{equation} \label{cosa}
\aligned b_{1,0} &= \int_0^t e^{\Delta(t-\tau)}\mathbb{P}\tilde
b_{1,0}(\cdot, \tau)d\tau \\
& = \frac {Q^2}{4r}e^{-t}\sin(\eta\cdot x)v\sum_{i=1}^r
|k_i||k_i'|\int_0^t e^{(1-(|k_i|^2 + |k_i'|^2)\tau} d\tau\\
& = \frac {Q^2}{4r}\sin(\eta\cdot x)v\sum_{i=1}^r |k_i||k_i'| \frac
{e^{-t} - e^{-(|k_i|^2 + |k_i'|^2)t}}{|k_i|^2 +
|k_i'|^2-1}\endaligned,\end{equation} 
since $v_i = v$ is fixed. Therefore we have
\begin{Lemma}\label{leb0}
Suppose $\frac{1}{|k_1|^2}\ll T\ll 1$. Then
\begin{equation}\label{b:b10}
\|b_{1,0}(\cdot, T)\|_{B_{\infty}^{-1, \infty}}=\sup_{\tau\in(0,
1)}\sqrt\tau \|e^{\tau\Delta}b_{1,0}\|_{L^\infty} \sim Q^2, \quad
\|b_{1,0}\|_{X_T} \lesssim\sqrt{T}Q^2.
\end{equation}
\end{Lemma}
\lpf By (\ref{cosa}), it follows that
\begin{equation} \label{b0-Q2}
b_{1,0} \sim Q^2\sin(\eta\cdot x)v,
\end{equation}
for $\frac{1}{|k_1|^2}\ll T\ll 1$. Indeed, $T\ll 1$ insures
$e^{-t}\sim 1$, for $t\leq T$; $\frac{1}{|k_1|^2}\ll T$ insures
$e^{-(|k_s|^2+|k'_s|^2)t}\sim 0$. And, since $|k_1|$ is very large,
$|k_s|\sim |k'_s|$ for every $s$ by (\ref{ks-ks'}). Thus,
\begin{align} \label{b1-Q2}
\|b_{1,0}\|_{B_{\infty}^{-1, \infty}}&\sim Q^2\sup_{0<t<1}\sqrt{t}\|e^{t\triangle}\sin(\eta\cdot x)\|_{L^{\infty}}\\
&\sim Q^2\sup_{0<t<1}\sqrt{t}e^{-|\eta|^2t}\notag\\
&\sim Q^2.\notag
\end{align}
And
\begin{align}
\|b_{1,0}\|_{X_T} &\sim Q^2\sup_{0<t<T}\sqrt{t}\|\sin(\eta\cdot x)\|_{L^{\infty}}\\
&+Q^2\sup_{x_0,0<R<T}\left(\frac{1}{|B(x_0,\sqrt{R})|}\int\limits_{0}^{R}
\int\limits_{B(x_0,\sqrt{R})}|\sin(\eta\cdot x)|^2dxdt\right)^{\frac{1}{2}}\notag\\
&\lesssim\sqrt{T}Q^2.\notag
\end{align}
\cbdu

\bigskip

\subsection{Analysis of $y$ and $z$}

In this section we analyze the parts $y$ and $z$ of the solution
The idea is to control $y$ and $z$ using the boundedness of the
bilinear operator $\mathcal{B}$ in the space $X_T$. Naively one
would hope that nonlinear terms turn out to be even smaller. But the
trouble is at the linear term $G_1$ and $K_1$ since
$$
\|e^{t\Delta}u_0\|_{X_T} \lesssim Q, \quad\text{and} \quad
\|e^{t\Delta}b_0\|_{X_T}\lesssim Q,
$$
by Lemma \ref{XTu0}, which is the best we can have. The problem is
the plane waves should not be lumped together in one single time scale.
Plane waves that have much bigger wave vectors diffuse much quicker.
Therefore it makes sense to analyze how $y$ and $z$ evolve in
different time scales and see how different plane waves contribute.
In \cite{BP}, Bourgain and Pavlovi\'{c} very skillfully designed 
time steps to group appropriately  the plane waves. We will use the same idea. We now introduce the time step division  as used in  \cite{BP}. Let
$$
0<T_1<T_2<\cdots <T_{\beta},
$$
where $\beta =Q^3$ and
$$
T_\alpha =|k_{r_\alpha}|^{-2},\;\; \;r_\alpha =r-\alpha Q^{-3}r, \ \
\ \alpha =1, 2, \ldots.
$$
In particular, $r_\beta =0$ and $T_\beta =|k_0|^{-2}$. The following
are the key estimates for  the time step design  in
\cite{BP}.
\begin{Lemma}\label{let} Suppose that $r$ is sufficiently large for
a fixed large number $Q$. Then
\begin{equation}\notag
\|(e^{t\triangle}u_0)\chi_{[T_\alpha, T_{\alpha
+1}]}(t)\|_{X_{T_{\alpha+1}}}\lesssim Q^{-1/2}, \quad
\|(e^{t\triangle}b_0)\chi_{[T_\alpha, T_{\alpha
+1}]}(t)\|_{X_{T_{\alpha+1}}}\lesssim Q^{-1/2}.
\end{equation}
\end{Lemma}
\lpf The proof is the same as  in \cite{BP}. For the
convenience of the reader we  outline the proof  showing 
 how the design of the time steps. First use the 
decomposition
\begin{equation}\notag
(e^{t\triangle}u_0)\chi_{[T_\alpha, T_{\alpha +1}]}(t)\approx
L_1+L_2+L_3,
\end{equation}
where
\begin{align}\notag
L_1&=\frac{Q}{\sqrt{r}}\sum_{s<r_{\alpha+1}}|k_s|v_s\cos(k_s\cdot x)e^{-|k_s|^2t}\chi_{[T_\alpha,T_{\alpha+1}]}(t)\\
L_2&=\frac{Q}{\sqrt{r}}\sum_{s=r_{\alpha+1}}^{r_\alpha}|k_s|v_s\cos(k_s\cdot x)e^{-|k_s|^2t}\chi_{[T_\alpha,T_{\alpha+1}]}(t)\notag\\
L_3&=\frac{Q}{\sqrt{r}}\sum_{r_\alpha<s\leq r}|k_s|v_s\cos(k_s\cdot
x)e^{-|k_s|^2t}\chi_{[T_\alpha,T_{\alpha+1}]}(t).\notag
\end{align}
The first group are those plane waves whose sizes are small.  Provided $\frac Q{\sqrt r} \lesssim Q^{-\frac 12}$ we have 

\begin{equation}\notag
\|L_1\|_{X_{T_{\alpha+1}}}
\lesssim\frac{Q}{\sqrt{r}}\sqrt{T_{\alpha+1}}|k_{r_{\alpha+1}-1}|+
\frac{Q}{\sqrt{r}}(T_{\alpha+1}|k_{r_{\alpha+1}-1}|^2)^{\frac{1}{2}}\notag
\lesssim\frac{Q}{\sqrt{r}} \lesssim Q^{-\frac 12},
\end{equation}
 The second group is
the group of active plane waves in the time scale $[T_\alpha,
T_{\alpha+1}]$. But, by the design, the number of plane waves in
this group is small.
\begin{align}\notag
\|L_2\|_{X_{T_{\alpha+1}}}
&\lesssim\frac{Q}{\sqrt{r}}+\frac{Q}{\sqrt{r}}(r_{\alpha}-r_{\alpha +1})^{\frac{1}{2}}\notag\\
&=\frac{Q}{\sqrt{r}}+\frac{Q}{\sqrt{r}}(Q^{-3}r)^{\frac{1}{2}}\notag\\
&\lesssim Q^{-1/2},\notag
\end{align}
if $\frac Q{\sqrt r} \lesssim Q^{-\frac 12}$. The last group of
plane waves are those that have diffused too much and become small
in size.
\begin{align}\notag
\|L_3\|_{X_{T_{\alpha+1}}} &\lesssim \frac{Q}{\sqrt{r}}\sup_{t<T_{\alpha+1}}\sum_{s=r_\alpha+1}^r|k_s| \sqrt{t}e^{-|k_s|^2t}\notag\\
&+\frac{Q}{\sqrt{r}}\sup_{t<T_{\alpha+1}}\left(\int\limits_0^t|\sum_{s=r_\alpha+1}^{r}|k_s|^2
e^{-|k_s|^2t}\chi_{[T_\alpha,T_{\alpha+1}]}(\tau)|d\tau\right)^{\frac{1}{2}}\notag
\end{align}
The first supremum of the last line is controlled by the integral
\begin{align}\notag
\int_{|k_{r_\alpha+1}|}^{|k_r|}\sqrt{t}e^{-x^2t}dx
&=\int_{|k_{r_\alpha+1}|\sqrt{t}}^{|k_r|\sqrt{t}}e^{-y^2}dy\\
&\leq\int_{|k_{r_\alpha+1}|\sqrt{t}}^{|k_r|\sqrt{t}}e^{-y}dy\notag\\
&\leq e^{-|k_{r_\alpha+1}|\sqrt{t}}\leq e^{-|k_{r_\alpha+1}|/|k_{r_\alpha}|}\ll 1,\notag
\end{align}
where we used the fact that $|k_{r_\alpha+1}|\sqrt{t}>1$ for $t\in[T_\alpha,T_{\alpha+1}]$ to the second step and the last step follows from (\ref{kl}).

The second supremum is controlled by
\begin{align}\notag
\left(\sum_{s=r_\alpha+1}^{r}e^{-|k_s|^2T_{\alpha+1}}\right)^{1/2}
&\leq (r-r_\alpha)e^{-|k_{r_\alpha+1}|^2/|k_{r_{\alpha+1}}|^2}\\
&\lesssim(\alpha Q^{-3}r)e^{-4^{r-\alpha Q^{-3}r}}\ll 1.\notag
\end{align}
In the same way one can show the same estimate for $b_0$. The proof
of Lemma is complete.\cbdu
\begin{Lemma}\label{T-beta} For $T>T_\beta$,
\begin{align}\label{u0XT2}
\|(e^{t\triangle}u_0)\chi_{[T_\beta,T]}(t)\|_{X_T}&\lesssim\frac{Q}{\sqrt{r}},
\end{align}
\begin{align}\label{b0XT2}
\|(e^{t\triangle}b_0)\chi_{[T_\beta,T]}(t)\|_{X_T}&\lesssim\frac{Q}{\sqrt{r}}.
\end{align}
\end{Lemma}
\lpf From Lemma \ref{k} and (\ref{h-u0}), we see that
\begin{align}
\|(e^{t\triangle}u_0)\chi_{[T_\beta,T]}(t)\|_{X_T}
&\lesssim\frac{Q}{\sqrt{r}}+ \frac{Q}{\sqrt{r}}
\sum_{s=1}^r|k_s|e^{-\frac{|k_s|^2}{|k_0|^2}}
|(T-T_\beta)^{\frac{1}{2}} \notag\\
&\lesssim\frac{Q}{\sqrt{r}}.\notag
\end{align}
The second one follows in the same way. \cbdu

Recall the equations for $y$ and $z$ from Section 2.4 that
\begin{align} \notag
y_t-\triangle y+G_0 + G_1 + G_2=0,\\
z_t-\triangle z+K_0 + K_1 + K_2=0,\notag
\end{align}
Note that, $y(0)=z(0)=0$. Hence, $t\in [T_\alpha, T_{\alpha+1}]$,
\begin{align}\label{y-t}
& y(t) = -\int\limits_{0}^{t}
e^{(t-\tau) \triangle} G (\tau)d\tau\\
& = -\int\limits_{0}^{t} e^{(t-\tau) \triangle}
G(\tau)\chi_{[0,T_\alpha]}(\tau)d\tau -\int\limits_{0}^{t}
e^{(t-\tau) \triangle} G(\tau)\chi_{[T_\alpha,
T_{\alpha+1}]}(\tau)d\tau, \notag
\end{align}
where $G = G_0 + G_1 + G_2$. So we can write
\begin{equation}\label{y}
\|y\|_{X_{T_{\alpha +1}}}\leq I_1+ I_2
\end{equation}
to see how $y$ develop in the time step $[T_\alpha, T_{\alpha+1}]$.

Similarly for $z$, we have
\begin{align}\label{z}
& \|z\|_{X_{T_{\alpha +1}}} \\
\leq
\|\int\limits_{0}^te^{(t-\tau)\triangle}K(\tau)\chi_{[0,T_\alpha]}d\tau\|_{X_{T_{\alpha
+1}}} & +
\|\int\limits_{0}^te^{(t-\tau)\triangle}K(\tau)\chi_{[T_\alpha,T_{\alpha+1}
]}d\tau\|_{X_{T_{\alpha +1}}} \notag \\
& = J_1+J_2,\notag
\end{align}
where $K = K_0 + K_1 + K_2$. Now we are ready to estimate the
increments of $y$ and $z$ during the time scale $[T_\alpha,T_{\alpha+1}]$.
\begin{Lemma}
With appropriate choice of $r$ and $T$, we have
\begin{align}\label{yzXT1}
\|y\|_{X_{T_{\alpha +1}}}+\|z\|_{X_{T_{\alpha +1}}}&\lesssim
Q^3(\frac{1}{r}+\sqrt{T_\beta})+Q(\|y\|_{X_{T_\alpha}}+\|z\|_{X_{T_\alpha}}).
\end{align}
\end{Lemma}
\lpf Applying bilinear estimate (\ref{B}), estimates in space $X_{T_\alpha}$ from Lemma \ref{XTu0},
(\ref{b:b10}), (\ref{b:b11}) and (\ref{XTu1}), we have
\begin{align}\label{I1}
I_1 &\lesssim
(Q+\frac{Q^2}{r}+\|y\|_{X_{T_\alpha}})\|y\|_{X_{T_\alpha}}\\
&+(Q+Q^2\sqrt{T_\alpha}+\frac{Q^2}{r}+\|z\|_{X_{T_\alpha}})\|z\|_{X_{T_\alpha}}\notag\\
&+(Q+\frac{Q^2}{r})\frac{Q^2}{r}\notag\\
&+(Q+Q^2\sqrt{T_\alpha}+\frac{Q^2}{r})(Q^2\sqrt{T_\alpha}+\frac{Q^2}{r}).\notag
\end{align}
We next apply the Lemma \ref{let} and estimate
\begin{align}\label{I2}
I_2
&\lesssim (Q^{-1/2}+\frac{Q^2}{r}+\|y\|_{X_{T_{\alpha+1}}})\|y\|_{X_{T_{\alpha+1}}}\\
&+(Q^{-1/2}+Q^2\sqrt{T_{\alpha+1}}+\frac{Q^2}{r}+\|z\|_{X_{T_{\alpha+1}}})\|z\|_{X_{T_{\alpha+1}}}\notag\\
&+(Q^{-1/2}+\frac{Q^2}{r})\frac{Q^2}{r}\notag\\
&+(Q^{-1/2}+Q^2\sqrt{T_{\alpha+1}}+\frac{Q^2}{r})(Q^2\sqrt{T_{\alpha+1}}+\frac{Q^2}{r}).\notag
\end{align}
We choose $r$ sufficiently large and $T$ appropriately small such that
\begin{equation}\label{qrt}
\frac{Q^2}{r}<Q^{-1/2}, \ \ Q^2\sqrt{T}<Q^{-1/2}.
\end{equation}
Hence we combine (\ref{y}), (\ref{I1}) and (\ref{I2}) to arrive at
\begin{eqnarray}\label{yXT1}
\|y\|_{X_{T_{\alpha +1}}}
\lesssim Q^3(\frac{1}{r}+\sqrt{T_\beta})+Q(\|y\|_{X_{T_\alpha}}+\|z\|_{X_{T_\alpha}})\\
+Q^{-1/2}(\|y\|_{X_{T_{\alpha+1}}}+\|z\|_{X_{T_{\alpha+1}}})+\|y\|_{X_{T_{\alpha
+1}}}^2+\|z\|_{X_{T_{\alpha +1}}}^2.\notag
\end{eqnarray}
Similarly we can obtain
\begin{align}\label{zXT1}
\|z\|_{X_{T_{\alpha +1}}} & \lesssim Q^3(\frac{1}{r}+\sqrt{T_\beta})+Q(\|y\|_{X_{T_\alpha}}+\|z\|_{X_{T_\alpha}})\\
&+Q^{-1/2}(\|y\|_{X_{T_{\alpha+1}}}+\|z\|_{X_{T_{\alpha+1}}})+\|y\|_{X_{T_{\alpha
+1}}}\|z\|_{X_{T_{\alpha +1}}}.\notag
\end{align}
Therefore, adding (\ref{yXT1}) and (\ref{zXT1}), we have
\begin{align}
\|y\|_{X_{T_{\alpha +1}}}& +\|z\|_{X_{T_{\alpha +1}}} \lesssim Q^3(\frac{1}{r}+\sqrt{T_\beta})
+Q(\|y\|_{X_{T_\alpha}}+\|z\|_{X_{T_\alpha}})\notag\\
& +(\|y\|_{X_{T_{\alpha +1}}}+\|z\|_{X_{T_{\alpha +1}}})^2.\notag
\end{align}
So, for much larger $r$ and $|k_0|$, we have $\|y\|_{X_{T_{\alpha
+1}}}+\|z\|_{X_{T_{\alpha +1}}}$ small and
$$
\|y\|_{X_{T_{\alpha +1}}}+\|z\|_{X_{T_{\alpha +1}}} \lesssim
Q^3(\frac{1}{r}+\sqrt{T_\beta})+Q(\|y\|_{X_{T_\alpha}}+\|z\|_{X_{T_\alpha}}).
$$ \cbdu

\bigskip

By iterating (\ref{yzXT1}) it then follows easily that
\begin{Lemma}
\begin{equation}\label{yzXT2}
\|y\|_{X_{T_\beta}}+\|z\|_{X_{T_\beta}}\lesssim Q^{\beta
+2}(\frac{1}{r}+\sqrt{T_\beta})
\end{equation}
\end{Lemma}
Next for $T>T_\beta$, in the light of Lemma \ref{T-beta}, one may
repeat the argument in the proof of (\ref{yzXT1}) and obtain

\begin{Lemma}\label{ley-xt}
For appropriate choice of $r$ and $T$,
\begin{equation}\label{y-xt}
\|y\|_{X_{T}} + \|z\|_{X_{T}} \lesssim  4Q^4T.
\end{equation}
\end{Lemma}
This implies that,
\begin{align}\label{yinf}
\|y(\cdot, T)\|_{B^{-1, \infty}_\infty} \lesssim \|y(\cdot,
T)\|_{L^\infty} \lesssim T^{-1/2}\|y\|_{X_T} \lesssim 4Q^4\sqrt{T}
\end{align}
and
\begin{align}\label{zinf}
\|z(\cdot, T)\|_{B^{-1, \infty}_\infty} \lesssim \|z(\cdot,
T)\|_{L^\infty} \lesssim T^{-1/2} \|z\|_{X_T} \lesssim 4Q^4\sqrt{T}.
\end{align}

\bigskip

\subsection{Finishing the Proof}

Now we are ready to complete the proof of Theorem \ref{theorem}.
Since (\ref{b:b11}) implies
\begin{equation}\label{b11inf}
\|b_{1,1}(\cdot, T)\|_{B^{-1, \infty}_\infty} \lesssim
\|b_{1,1}(\cdot, T)\|_{L^\infty}\lesssim T^{-\frac{1}{2}}\|b_{1,1}\|_{X_T} \lesssim \frac{Q^2}{r\sqrt{T}},
\end{equation}
from (\ref{b}) we combine (\ref{b:b10}), (\ref{b11inf}) and
(\ref{zinf}) to obtain that
\begin{align}
\|b(\cdot, T)-e^{T\triangle}b_0\|_{\dot{B}_\infty^{-1,\infty}} & \geq
\|b(\cdot, T)-e^{T\triangle}b_0\|_{B_\infty^{-1,\infty}} \notag\\
&\geq\|b_{1, 0}(\cdot, T)\|_{B_\infty^{-1,\infty}} - \|b_{1,1}(\cdot, T)\|_{B_\infty^{-1,\infty}}
- \|z(\cdot, T)\|_{B_\infty^{-1,\infty}} \notag \\
& \gtrsim Q^2-\|b_{1,1}\|_{L^\infty}-\|z\|_{L^\infty}\notag\\
& \gtrsim Q^2(1-\frac{1}{r\sqrt{T}}-4Q^2\sqrt{T}).\notag
\end{align}
Therefore, in the light of (\ref{norm:u0-3}),
\begin{equation}\notag
\|b(T)\|_{\dot{B}_\infty^{-1,\infty}} \gtrsim Q^2.
\end{equation}
On the other hand, from (\ref{u}), we combine (\ref{norm:u0-3}),
(\ref{XTu1}), and (\ref{yinf}) and have, for any $t\in [0, T]$,
\begin{equation}
\|u(\cdot, t)\|_{B^{-1, \infty}_\infty} \lesssim \frac Q{\sqrt r} +
\frac {Q^2}r + Q^4\sqrt T
\end{equation}
and remains small in $B^{-1, \infty}_\infty$. Thus we proved theorem
\ref{theorem}.

\begin{Remark}
We would like to note the following simple chart to indicate how
the choices of the several parameters are made.
\begin{equation}\notag
\delta\longrightarrow Q\longrightarrow T\longrightarrow
|k_0|\longrightarrow |k_s|, \ \ Q\longrightarrow r.
\end{equation}
\end{Remark}

\bigskip

\section{Other scenarios of norm inflations}

In this section, we consider other interesting norm inflation
phenomena for MHD system. The essence of the construction introduced
by Bourgain and Pavlovi\'{c} in \cite{BP} and used in this paper is
that the collisions of plane waves with similar but large wave
vectors can cause norm inflations in NSE as well as in MHD system.
More precisely we see that the collisions in the quadratic terms
trigger the norm inflations. Hence we can arrange the initial data
to have collisions of plane waves with similar wave vectors either
in $u_1$ or in $b_1$ to produce various norm inflation modes for
MHD system. For example, even the initial velocity is zero, if the
initial magnetic field contains enough plane waves to collide, 
we can produce the scenario where the velocity develops norm
inflation in $\dot{B}^{-1, \infty}_\infty$ while the magnetic field remains small in the space
$B^{-1, \infty}_\infty$. Namely,
\begin{Theorem}\label{theorem u0-0}
Let $u_0\equiv 0$ and
$$
b_0 = \frac Q{\sqrt r}\sum_{i=1}^r (|k_i|v_i\cos(k_i\cdot x) +
|k_i'|v_i'\cos(k_i'\cdot x)).
$$
Then, for any $\delta>0$, there exists a solution $(u,b,p)$ to the
MHD system (\ref{MHD:PDE1}) with initial data $u_0$ and $b_0$ as in
the above satisfying
\begin{equation}\notag
\|b(0)\|_{\dot{B}_{\infty}^{-1, \infty}} \lesssim \delta,
\end{equation}
such that for some $0<T<\delta$
\begin{equation}\notag
\|u(T)\|_{\dot{B}_{\infty}^{-1, \infty}} \gtrsim 1/ \delta,
\end{equation}
while for any $0<t<T<\delta$
\begin{equation}\notag
\|b(t)\|_{{B}_{\infty}^{-1, \infty}} \lesssim \delta.
\end{equation}
\end{Theorem}
The proof will be more or less the same as the proof in \cite{BP} in
the light of our discussions in the previous section. Another
interesting case mentioned in Remark \ref{remark2} is that when the
initial velocity and the initial magnetic field are the same, although
they both include many plane waves that are to collide, the
collisions cancel each other in the evolution of the MHD system and
produce no norm inflations.
\begin{Remark}
Finally we would like to mention that, in \cite{CS},  Cheskidov and
Shvydkoy introduced a different construction of initial data to
prove the ill-posedness of NSE in certain Besov spaces. There are two more works
about the ill-posedness results for Navier-Stokes equations by Germain \cite{Ger}
and Yoneda \cite{Yo}.
\end{Remark}

\bigskip

\noindent {\bf Acknowledgements} The authors would like to thank
Natasa Pavlovic for many helpful remarks and comments and 
Lorenzo Brandolese for some useful suggestions.

\end{document}